# Orbits of Quaternionic Möbius Transformations

# Tony Thrall

# Abstract

Möbius transformations of the extended complex plane are at the crossroads of many interesting topics, e.g., they form a group under composition, are the simplest form of rational function, and are a path to Lie theory. Quaternionic transformations are a subgroup of Möbius transformations isomorphic to rotations of the Riemann sphere, which also represent quaternion conjugation. These representations yield formulas for the axis and radians of rotation, and thereby portray each particular quaternionic transformation as part of a continuous orbit of rotations sharing a common axis.

# Introduction

The purpose of this note is to show that each quaternionic Möbius transformation is the cumulative change induced by a continuous motion. We first provide a quick tour of the components that we will assemble in the next section. These components are presented more comprehensively in the references cited below: for Möbius transformations and the Riemann sphere see chapter 3 of [3]; for quaternions see chapters 3 and 4 of [2], along with chapter 1 of [5]; and for relevant Lie theory see [1], [4], and [5].

## Möbius transformations

A Möbius transformation $M(\cdot) = \mathcal{M}(\cdot, a, b, c, d)$ has the following form:

$$M(z) = \frac{a\,z+b}{c\,z+d} \qquad \text{with } z, a, b, c, d \in \mathbb{C} \text{ and with } a\,d - b\,c \neq 0$$

We adopt the following definitions to expand the domain and range of $M(\cdot)$ to the extended complex plane, $\hat{\mathbb{C}}$, which is the union of $\mathbb{C}$ and the point at infinity: $\hat{\mathbb{C}} \overset{def}{=} \mathbb{C} \cup \{\infty\}$.

$$M(\infty) \overset{def}{=} \begin{cases} a/c & \text{if } c \neq 0 \\ \infty & \text{if } c = 0 \end{cases}$$

$$M(-d/c) \overset{def}{=} \begin{cases} \infty & \text{if } c \neq 0 \\ M(\infty) & \text{if } c = 0 \end{cases}$$

Each point $z \in \hat{\mathbb{C}}$ can be represented in homogeneous coordinates as the equivalence class of pairs of complex values, $\langle z, 1\rangle = \{k\,z, k\}_k$ for finite $z$ or $\langle 1, 0\rangle = \{k, 0\}_k$ for $z = \infty$, in both cases with complex $k \neq 0$.

In this setting $M$ sends the $\langle z, 1\rangle$ equivalence class to the $\langle a z + b, c z + d\rangle$ equivalence class, and is thus represented by the linear transformation

$$\left\langle\begin{pmatrix} z \\ 1 \end{pmatrix}\right\rangle \mapsto \left\langle\begin{pmatrix} a & b \\ c & d \end{pmatrix}\begin{pmatrix} z \\ 1 \end{pmatrix}\right\rangle$$

Then the composition of transformations $M_2(M_1(z))$ can be shown to correspond to the composition of linear operators on homogeneous coordinates, and thus to matrix multiplication with respect to homogeneous coordinates. Continuing on this path one can show that Möbius transformations are a group under the composition operation.

Note that for complex $k \neq 0$, $\mathcal{M}(\cdot, a, b, c, d)$ and $\mathcal{M}(\cdot, k a, k b, k c, k d)$ refer to the same transformation. This ambiguity is reduced if we normalize the coefficients of $M$ so that $a d - b c = 1$, but even then $\mathcal{M}(\cdot, a, b, c, d)$ and $\mathcal{M}(\cdot, -a, -b, -c, -d)$ refer to the same transformation. (Formally, the restriction of $\mathcal{M}$ to normalized coefficients can be represented as an isomorphism from the quotient group $SL(2, \mathbb{C})/\{\pm 1\}$ of the special linear group to the group of Möbius transformations.)

With normalized coefficients the inverse of $M$ is $M^{-1}(\cdot) = \mathcal{M}(\cdot, d, -b, -c, a)$. It is also useful to borrow from matrices the conjugate-transpose notation by defining the transformation-conjugate of $M$ as $M^*(\cdot) = \mathcal{M}(\cdot, \overline{a}, \overline{c}, \overline{b}, \overline{d})$.

# The Riemann sphere

The two-dimensional unit sphere, $S^2 = \{\boldsymbol{\eta} \in \mathbb{R}^3 \mid \|\boldsymbol{\eta}\| = 1\}$, is known as the Riemann sphere when mapped bijectively to the extended complex plane $\hat{\mathbb{C}}$ via stereographic projection, $\sigma$.

$$\sigma(\eta_1, \eta_2, \eta_3) = \frac{\eta_1 + i\,\eta_2}{1 - \eta_3} \quad \text{with } \sigma(0, 0, 1) = \infty$$

$$\sigma^{-1}(x + i\,y) = \left\{\frac{2x}{x^2+y^2+1}, \frac{2y}{x^2+y^2+1}, \frac{x^2+y^2-1}{x^2+y^2+1}\right\} \text{with } \sigma^{-1}(\infty) = \{0, 0, 1\}$$

This mapping of the extended plane to a compact surface is useful for visualizing Möbius transformations (and other types of complex-valued functions of a complex variable) as amply displayed in [3].

# Quaternions

Quaternions, denoted by $\mathbb{H}$ in honor of the Irish mathematician W. R. Hamilton (1805-1865), are a real four-dimensional algebra, whose standard basis we will denote as $\{1, i_1, i_2, i_3\}$, so that any quaternion $h \in \mathbb{H}$ has a unique set of coordinates with respect to this basis.

$$h = h_0 + h_1 i_1 + h_2 i_2 + h_3 i_3 \qquad \text{with } h_0, h_1, h_2, h_3 \in \mathbb{R}$$

We refer to $h_0$ as the real part of $h$, denoted $\mathcal{R}e(h)$, and to $(h_1 i_1 + h_2 i_2 + h_3 i_3)$ as the vector part of $h$, denoted $\mathcal{V}ec(h)$. The related quaternion, $h_0 - h_1 i_1 - h_2 i_2 - h_3 i_3$, is known as the quaternion conjugate of $h$, and is denoted as $h^*$. That is,

$$h = \mathcal{R}e(h) + \mathcal{V}ec(h)$$

$$h^* = \mathcal{R}e(h) - \mathcal{V}ec(h)$$

Multiplication in $\mathbb{H}$ is in general not commutative, being determined by the following equalities.

$$-1 = i_1^2 = i_2^2 = i_3^2 = i_1\, i_2\, i_3$$

$\mathbb{H}$ is also equipped with the inner product and norm of $\mathbb{R}^4$, i.e., for $h, q \in \mathbb{H}$ we have

$$<h, q> = h_0\, q_0 + h_1\, q_1 + h_2\, q_2 + h_3\, q_3 \quad \text{and} \quad \|h\|^2 = <h, h>$$

The inner product may be defined without reference to the standard coordinates as $<h, q> = \mathcal{R}e(h\, q^*)$. A quaternion $v \in \mathbb{H}$ whose real part is equal to zero, $\mathcal{R}e(v) = 0$, is said to be a pure vector quaternion, collectively denoted as $\mathcal{V}ec(\mathbb{H})$. We denote by $\mathcal{U}(\mathbb{H})$ the subset of pure vectors $u$ having unit norm. It turns out that a quaternion $u$ is member of $\mathcal{U}(\mathbb{H})$ if and only if $u^2 = -1$. The members of $\mathcal{U}(\mathbb{H})$ play the role of the imaginary unit, $i \in \mathbb{C}$. For example, just as non-zero $z \in \mathbb{C}$ can be expressed in polar coordinates, so can non-zero $q \in \mathbb{H}$.

$$q = \|q\|\, e^{u\,\theta} = \|q\|\,(\cos(\theta) + u \sin(\theta))$$

for certain $\theta \in \mathbb{R}$ and $u \in \mathcal{U}(\mathbb{H})$, namely

$$\theta = \arctan(\mathcal{R}e(q),\ \|\mathcal{V}ec(q)\|)$$

$$u = \mathcal{V}ec(q) / \|\mathcal{V}ec(q)\| \quad \text{if } \mathcal{V}ec(q) \neq 0$$

If $\mathcal{V}ec(q) = 0$, then $u$ may be any member of $\mathcal{U}(\mathbb{H})$. The polar representation of non-zero $q$ clarifies the effect of quaternion conjugation $\Gamma_q(\cdot)$:

$$\Gamma_q : h \mapsto q\, h\, q^{-1} = e^{u\,\theta}\, h\, e^{-u\,\theta} \qquad \text{for } h, q \in \mathbb{H} \text{ with } q \neq 0$$

One can determine $v, w \in \mathcal{U}(\mathbb{H})$ so that $\{u, v, w\}$ is an orthonormal basis for $\mathcal{V}ec(\mathbb{H})$ with the same multiplication rules as $\{i_1, i_2, i_3\}$, namely

$$-1 = u^2 = v^2 = w^2 = u\, v\, w$$

If we now express $h$ in the $\{u, v, w\}$-coordinate system

$$h = h_0 + h_u\, u + h_v\, v + h_w\, w$$

then one can show that

$$h \mapsto q\, h\, q^{-1} = h_0 + h_u\, u + h_v\, v_\theta + h_w\, w_\theta$$

where

$$\begin{pmatrix} v_\theta \\ w_\theta \end{pmatrix} = \begin{pmatrix} \cos(2\,\theta) & \sin(2\,\theta) \\ -\sin(2\,\theta) & \cos(2\,\theta) \end{pmatrix} \begin{pmatrix} v \\ w \end{pmatrix}$$

That is, conjugation of $h$ by $q$ rotates $h$ about axis $u$ by $2\,\theta$ radians, as presented more fully in [2] and [5]. In the next section we note that the subject of our discussion, quaternionic Möbius transformations, rotate the Riemann sphere, a fact that enables us to identify these transformations with quaternion conjugation.

# Quaternionic Möbius transformations

In the remainder of this note we focus on quaternionic transformations, which are of the form $M(\cdot) = \mathcal{M}(\cdot, \zeta, -\omega, \overline{\omega}, \overline{\zeta})$, and which we will abbreviate as $Q(\cdot, \zeta, \omega)$.

$$M(\cdot) = Q(\cdot, \zeta, \omega) \overset{def}{=} \mathcal{M}(\cdot, \zeta, -\omega, \overline{\omega}, \overline{\zeta})$$

Thus

$$M(z) = \frac{\zeta z - \omega}{\overline{\omega} z + \overline{\zeta}}$$

We will assume the coefficients are normalized, that is, that they satisfy $|\zeta|^2 + |\omega|^2 = 1$. Then

$$M^{-1} = M^* = Q(\cdot, \overline{\zeta}, -\omega) \text{ for } |\zeta|^2 + |\omega|^2 = 1$$

and $Q$ can be represented as an isomorphism from the quotient group SU(2, $\mathbb{C}$)/{±1} of the special unitary group to the group of quaternionic Möbius transformations. The term quaternionic derives from the fact that the matrices appearing in the transformation of homogeneous coordinates

$$\left\langle\!\left\langle \begin{pmatrix} z \\ 1 \end{pmatrix} \right\rangle\!\right\rangle \mapsto \left\langle\!\left\langle \begin{pmatrix} \zeta & -\omega \\ \overline{\omega} & \overline{\zeta} \end{pmatrix} \begin{pmatrix} z \\ 1 \end{pmatrix} \right\rangle\!\right\rangle$$

along with matrix multiplication, are isomorphic to the multiplicative group of quaternions having unit norm. Transformations of this form, under inverse stereographic projection to the Riemann sphere, can be shown to map each pair of antipodal points to another pair of antipodal points ([3]). Since the transformations are continuous, one can then conclude that they correspond to rotations of the Riemann sphere. The matrix $[\hat{M}]$ of this rotation of $\mathbb{R}^3$, in standard coordinates, is as follows.

$$[\hat{M}] = \begin{pmatrix} \mathcal{Re}(\zeta^2 - \omega^2) & -\mathcal{Im}(\zeta^2 + \omega^2) & 2\,\mathcal{Re}(\zeta\,\omega) \\ \mathcal{Im}(\zeta^2 - \omega^2) & \mathcal{Re}(\zeta^2 + \omega^2) & 2\,\mathcal{Im}(\zeta\,\omega) \\ -2\,\mathcal{Re}(\zeta\,\overline{\omega}) & 2\,\mathcal{Im}(\zeta\,\overline{\omega}) & |\zeta|^2 - |\omega|^2 \end{pmatrix}$$

Quaternion conjugation $\Gamma_q : h \mapsto q\,h\,q^{-1}$, which rotates $\mathcal{Vec}(h)$, can also be expressed as a rotation of 3-space. If we assume $q$ to have norm 1, then the matrix $[C_q]$ of this rotation in standard coordinates is as follows.

$$[C_q] = \begin{pmatrix} q_0^2 + q_1^2 - q_2^2 - q_3^2 & 2(q_1 q_2 - q_0 q_3) & 2(q_0 q_2 + q_1 q_3) \\ 2(q_0 q_3 + q_1 q_2) & q_0^2 - q_1^2 + q_2^2 - q_3^2 & 2(q_2 q_3 - q_0 q_1) \\ 2(q_1 q_3 - q_0 q_2) & 2(q_0 q_1 + q_2 q_3) & q_0^2 - q_1^2 - q_2^2 + q_3^2 \end{pmatrix} \quad \text{where } \|q\| = 1$$

If we identify these two rotations we obtain the following bijective maps $\gamma$ and $\gamma^{-1}$ between quaternions $q$ and the parameters $\{\zeta, \omega\}$ of transformations $M$.

$$\{\zeta, \omega\} = \{\zeta_q, \omega_q\} = \gamma(q) \overset{def}{=} \{q_0 + i\, q_3,\ q_2 - i\, q_1\}$$

$$q = q_{\zeta,\omega} = \gamma^{-1}\{\zeta, \omega\} = \mathcal{Re}(\zeta) - \mathcal{Im}(\omega)\, i_1 + \mathcal{Re}(\omega)\, i_2 + \mathcal{Im}(\zeta)\, i_3$$

In fact $\gamma$ yields a group isomorphism $\tilde{\gamma}$ from: (1) rotations of $\mathcal{Vec}(\mathbb{H})$ induced by quaternion conjugation $\Gamma_q$; to (2) quaternionic Möbius transformations $Q(\cdot, \gamma(q)) = Q(\cdot, \zeta, \omega)$.

# Polar decomposition

The axis of quaternion rotation can be defined as the line through $\pm u$, whose standard vector coordinates we denote as $[\pm u]$. These are antipodal points of the Riemann sphere, which we identify as the axis of rotation induced by $M$. The stereographic projection of these antipodal points, $\sigma([\pm u])$, equals the set of fixed points of $M$.

We will now decompose this rotation of the Riemann sphere as: (1) a rotation of $[u]$ to the north pole {0, 0, 1}; (2) a rotation about the polar axis through $2\,\theta$ radians; and (3) a return rotation of the north pole to $[u]$. This decomposition, presented in [4] for rotations about the coordinate axes, amounts to a conjugation of rotations, with the central rotation about the north pole representing the complex variable transformation $z \mapsto e^{2 i \theta} z$. We will use this rotational decomposition to obtain a polar decomposition of transformation $M$.

## A quaternion basis aligned with the given rotation

The decomposition outlined above begins with the rotation of the Riemann sphere that takes $[u]$ to the north pole $\{0, 0, 1\}$. The axis of rotation must therefore be orthogonal to both $[u]$ and $\{0, 0, 1\}$, and can thus be defined as the normalized cross-product of these two vectors. Let $w$ denote the unit vector quaternion whose standard coordinates $[w]$ are equal to those of the normalized cross-product.

$$[w] = \frac{[u]\otimes\{0,0,1\}}{\|[u]\otimes\{0,0,1\}\|} = \frac{\{u_2, -u_1, 0\}}{\sqrt{u_1^2+u_2^2}}$$

We can now construct $v \in \mathcal{U}(\mathbb{H})$ so that $\{u, v, w\}$ constitute an orthonormal basis for $\mathcal{V}ec(\mathbb{H})$, as previously described. That is, we define $[v]$ to be the cross-product $[w]\otimes[u]$.

$$[v] = [w]\otimes[u] = \frac{\{-u_1 u_3, -u_2 u_3, u_1^2+u_2^2\}}{\sqrt{u_1^2+u_2^2}}$$

## The given rotation as a conjugation of rotations

The angle of the initial component rotation can be determined by expressing $[u]$ in terms of its spherical coordinates, say $\{\sin(\phi)\cos(\chi), \sin(\phi)\sin(\chi), \cos(\phi)\} = \{u_1, u_2, u_3\}$, which is to be rotated to the north pole, with coordinates $\{0, 0, \cos(0)\}$. Thus the angle $\phi$ is

$$\phi = \arctan\left(u_3, \sqrt{u_1^2 + u_2^2}\right)$$

and the angle of rotation that takes $[u]$ to the north pole is $-\phi$. Consequently this rotation can be expressed as the following quaternion conjugation about axis $w$ through $-\phi$ radians.

$$h \mapsto e^{-w\,\phi/2}\, h\, e^{w\,\phi/2}$$

The second rotation is about the polar axis through an angle of $2\theta$ radians. Since the north pole has $[i_3]$ as its standard coordinates, the second spherical rotation represents the following quaternion rotation.

$$h \mapsto e^{i_3\,\theta}\, h\, e^{-i_3\,\theta}$$

The final rotation is the inverse of the first rotation, i.e.,

$$h \mapsto e^{w\,\phi/2}\, h\, e^{-w\,\phi/2}$$

## Corresponding polar decomposition

We now use the mapping $\gamma$ and its inverse to translate the preceding composition of three quaternion rotations into the composition of three quaternionic Möbius transformations. We begin by expressing the axis of rotation, $u$, as a function $u_{\zeta,\omega}$ of parameters $\{\zeta, \omega\}$ via $\gamma^{-1}$. Recalling that $|\zeta|^2 + |\omega|^2 = 1$, we have

$$u = \mathcal{V}ec(q) / \|\mathcal{V}ec(q)\| = \frac{-\mathcal{I}m(\omega)\, i_1 + \mathcal{R}e(\omega)\, i_2 + \mathcal{I}m(\zeta)\, i_3}{\|-\mathcal{I}m(\omega)\, i_1 + \mathcal{R}e(\omega)\, i_2 + \mathcal{I}m(\zeta)\, i_3\|}$$

$$= \frac{-\mathcal{I}m(\omega)\, i_1 + \mathcal{R}e(\omega)\, i_2 + \mathcal{I}m(\zeta)\, i_3}{\sqrt{|\omega|^2 + \mathcal{I}m(\zeta)^2}} = \frac{-\mathcal{I}m(\omega)\, i_1 + \mathcal{R}e(\omega)\, i_2 + \mathcal{I}m(\zeta)\, i_3}{\sqrt{1 - \mathcal{R}e(\zeta)^2}}$$

$$\overset{def}{=} u_{\zeta,\omega}$$

Similarly,

$$\theta = \arctan(\mathcal{R}e(q),\ \|\mathcal{V}ec(q)\|) = \arctan\left(\mathcal{R}e(\zeta),\ \sqrt{1 - \mathcal{R}e(\zeta)^2}\right)$$

from which we define $\tau_\zeta = 2\,\theta \bmod 2\,\pi$ as the radians of rotation induced by $M$.

$$\tau_\zeta \overset{def}{=} 2\arctan\left(\mathcal{R}e(\zeta),\ \sqrt{1 - \mathcal{R}e(\zeta)^2}\right) \bmod 2\,\pi \quad \text{for } |\zeta| \le 1$$

Continuing with this development, we define $w_\omega$, $\phi_\zeta$, and $v_{\zeta,\omega}$ as follows.

$$w = \frac{u_2\, i_1 - u_1\, i_2}{\sqrt{u_1^2 + u_2^2}} = \frac{\mathcal{R}e(\omega)\, i_1 + \mathcal{I}m(\omega)\, i_2}{|\omega|}$$

$$\overset{def}{=} w_\omega \qquad \text{for } 0 < |\omega| \le 1$$

$$\phi = \arctan\left(u_3,\ \sqrt{u_1^2 + u_2^2}\right)$$

$$= \arctan(\mathcal{I}m(\zeta),\ |\omega|) = \arctan\left(\mathcal{I}m(\zeta),\ \sqrt{1 - |\zeta|^2}\right)$$

$$\overset{def}{=} \phi_\zeta \qquad \text{for } |\zeta| \le 1 \text{ and } \zeta \ne \pm 1$$

$$v = \frac{\{-u_1\, u_3, -u_2\, u_3, u_1^2 + u_2^2\}}{\sqrt{u_1^2 + u_2^2}} = \frac{\mathcal{I}m(\zeta)\, \mathcal{I}m(\omega)\, i_1 - \mathcal{I}m(\zeta)\, \mathcal{R}e(\omega)\, i_2 + |\omega|^2\, i_3}{|\omega|\, \sqrt{1 - \mathcal{R}e(\zeta)^2}}$$

$$\overset{def}{=} v_{\zeta,\omega} \text{ for } \omega \ne 0 \text{ and } |\zeta|^2 + |\omega|^2 = 1$$

Then, for $\omega \ne 0$ and $|\zeta|^2 + |\omega|^2 = 1$, we have

$$\zeta = \cos(\tau_\zeta/2) + i \sin(\tau_\zeta/2) \cos(\phi_\zeta)$$

$$\omega = \sin(\tau_\zeta/2) \sin(\phi_\zeta)\, e^{i \arg(\omega)}$$

and the quaternion conjugation

$$h \mapsto e^{u_{\zeta,\omega}\, \tau_\zeta/2}\, h\, e^{-u_{\zeta,\omega}\, \tau_\zeta/2}$$

corresponds to $M = Q(\cdot, \zeta, \omega)$. To identify the components of the composition of transformations, we begin with the third component quaternion conjugation

$$h \mapsto e^{w_\omega\, \phi_\zeta/2}\, h\, e^{-w_\omega\, \phi_\zeta/2}$$

and we will denote its corresponding Möbius transformation as $W(\cdot)$, whose parameters are obtained from $\gamma$.

$$e^{w_\omega\, \phi_\zeta/2} = \cos(\phi_\zeta/2) + \sin(\phi_\zeta/2)\, w_\omega$$

$$= \cos(\phi_\zeta/2) + \sin(\phi_\zeta/2)\, (\mathcal{R}e(\omega)\, i_1 + \mathcal{I}m(\omega)\, i_2) / |\omega|$$

$$\gamma\left(e^{w_\omega\, \phi_\zeta/2}\right) = \{\cos(\phi_\zeta/2),\ \sin(\phi_\zeta/2)\, (\mathcal{I}m(\omega) - i\, \mathcal{R}e(\omega)) / |\omega|\}$$

$$= \left\{\cos(\phi_\zeta/2),\ \sin(\phi_\zeta/2)\, \frac{-i\, \omega}{|\omega|}\right\}$$

$$W(\cdot) \overset{def}{=} Q\left(\cdot, \cos(\phi_\zeta/2), -\sin(\phi_\zeta/2)\, \frac{i\,\omega}{|\omega|}\right)$$

Then the first component conjugation

$$h \mapsto e^{-w_\omega \phi_\zeta/2}\, h\, e^{w_\omega \phi_\zeta/2}$$

corresponds to

$$W^*(\cdot) = Q\left(\cdot, \cos(\phi_\zeta/2), \sin(\phi_\zeta/2)\, \frac{i\,\omega}{|\omega|}\right)$$

Finally, as previously noted, the central component conjugation

$$h \mapsto e^{i_3 \tau_\zeta/2}\, h\, e^{-i_3 \tau_\zeta/2}$$

corresponds to the complex variable transformation $z \mapsto e^{i\,\tau_\zeta} z$, which we denote as the quaternion Möbius transformation $D$ (corresponding to a diagonal matrix).

$$D(\cdot) \overset{def}{=} Q\left(\cdot, e^{i\,\tau_\zeta/2}, 0\right)$$

Based on the isomorphism provided by $\gamma$, we have the following decomposition of $M$.

$$M = W \circ D \circ W^*$$

We refer to this as a polar decomposition since in the matrix version of this identity

$$[M] = [W]\,[D]\,[W]^*$$

the unitary matrix $[M]$ is expressed as the conjugation of a diagonal matrix $[D]$ by a unitary matrix $[W]$.

To recap, the transformation $W$ corresponds to a spherical rotation that takes the north pole to $[u_{\zeta,\omega}]$, which defines the axis of $M$'s induced rotation of the Riemann sphere. The number of radians of $M$'s rotation is encoded in the transformation $D$. The separation of these two defining aspects of $M$'s rotation enables us, in the next section, to identify a continuous orbit that contains $M$ along with other quaternion Möbius transformations.

# A continuum of transformations

## Enveloping single-parameter sets

The polar decomposition enables us to portray the given transformation $M$ as one point in an orbit $\{\mathcal{T}_\tau\}$, as follows.

$$\mathfrak{D}^{\circ\tau}(\cdot) \overset{def}{=} Q\left(\cdot, e^{i\,\tau/2}, 0\right) : z \mapsto e^{i\,\tau} z$$

$$\left\{\mathcal{T}_\tau \overset{def}{=} W \circ \mathfrak{D}^{\circ\tau} \circ W^* \;\middle|\; \tau \in \mathbb{R}\right\} \quad \text{with } D = \mathfrak{D}^{\circ\tau_\zeta} \text{ and } M = \mathcal{T}_{\tau_\zeta}$$

(The notation $\mathfrak{D}^{\circ\tau}$ generalizes the $n$-fold composition $\mathfrak{D}^{\circ n}$ of $\mathfrak{D}^{\circ 1}$ with itself $n$ times, for integer $n$.)

$\mathcal{T}_\tau$ represents a spherical rotation of $\tau$ radians about the same axis of rotation $[u_{\zeta,\omega}]$ as $M$. In the extended complex plane, $\mathcal{T}_\tau$ has the same fixed points and invariant curves as $M$, so that for each $z \in \hat{\mathbb{C}}$, $\{\mathcal{T}_\tau(z)\}$ is the invariant curve containing $z$. $\{\mathcal{T}_\tau\}$ is a commutative subgroup of quaternionic Möbius transformations, with $\mathcal{T}_\tau \circ \mathcal{T}_\upsilon = \mathcal{T}_{\tau+\upsilon}$, $\mathcal{T}_{\tau+2\pi} = \mathcal{T}_\tau$, and $\mathcal{T}_0$ the identity transformation.

We have shown $M$ to be a member of $\{\mathcal{T}_\tau\}$, in which we fixed the axis of rotation $u_{\zeta,\omega}$ and varied the radians of rotation $\tau$. We can just as well fix the radians of rotation and vary the axis of rotation. If we do so along the great circle from the north pole through $[u_{\zeta,\omega}]$ we obtain

$$\left\{\Phi_\phi \overset{def}{=} \mathfrak{W}^{\circ\phi} \circ D \circ \mathfrak{W}^{\circ -\phi} \;\middle|\; \phi \in \mathbb{R}\right\} \text{ with } W = \mathfrak{W}^{\circ\phi_\zeta} \text{ and } M = \Phi_{\phi_\zeta}$$

where $\mathfrak{W}^{\circ\phi}$ denotes the transformation corresponding to the quaternion conjugation $h \mapsto e^{w\,\phi/2}\, h\, e^{-w\,\phi/2}$.

$$\mathfrak{W}^{\circ\phi}(\cdot) \overset{def}{=} Q\left(\cdot,\ \cos\left(\frac{\phi}{2}\right),\ -\sin\left(\frac{\phi}{2}\right) \frac{i\,\omega}{|\omega|}\right)$$

The set of transformations $\{\mathfrak{W}^{\circ\phi}\}$ is a commutative subgroup, with $\mathfrak{W}^{\circ\psi} \circ \mathfrak{W}^{\circ\phi} = \mathfrak{W}^{\circ(\psi+\phi)}$ and $\mathfrak{W}^{\circ(\phi+2\pi)} = \mathfrak{W}^{\circ\phi}$, since $Q(\alpha, \beta)$ and $Q(-\alpha, -\beta)$ refer to the same transformation. The set $\{\Phi_\phi\}$ is generally not a subgroup, however, unless $D$ is the identity transformation, in which case $\{\Phi_\phi\}$ consists solely of the identity transformation. Parameter $\phi$ gives the angular declination from the north pole of the Riemann sphere along the great circle through $[u_{\zeta,\omega}]$.

A third set containing $M$ retains both the degree of rotation, $\tau = \tau_\zeta$, and the declination, $\phi = \phi_\zeta$, of the axis of rotation but varies the longitude $\lambda$ of the axis of rotation along the fixed latitude defined by $\phi$. We therefore define

$$\mathfrak{L}_\lambda(\cdot) \overset{def}{=} Q\left(\cdot,\ \cos\left(\frac{\phi_\zeta}{2}\right),\ -\sin\left(\frac{\phi_\zeta}{2}\right) e^{i\,\lambda}\right) \quad \text{with } W = \mathfrak{L}_{\arg(i\,\omega)}$$

$$\left\{\Lambda_\lambda \overset{def}{=} \mathfrak{L}_\lambda \circ D \circ \mathfrak{L}_\lambda^* \;\middle|\; \lambda \in \mathbb{R}\right\} \qquad \text{with } M = \Lambda_{\arg(i\,\omega)}$$

This third set $\{\Lambda_\lambda\}$ has $2\pi$ periodicity and each pair of elements commute, i.e., $\Lambda_\mu \circ \Lambda_\lambda = \Lambda_\lambda \circ \Lambda_\mu$, but this product is in general not a member of $\{\Lambda_\lambda\}$. That is, $\{\Lambda_\lambda\}$ is in general not a subgroup, unless $D$ is the identity transformation, in which case $\{\Lambda_\lambda\}$ consists solely of the identity transformation.

Of course any axis of rotation can be determined from a combination of the declination parameter $\phi$ and the longitude parameter $\lambda$, which leads to the following family of transformations.

$$\mathfrak{U}_{\phi,\lambda}(\cdot) \overset{def}{=} Q\left(\cdot,\ \cos\left(\frac{\phi}{2}\right),\ -\sin\left(\frac{\phi}{2}\right) e^{i\,\lambda}\right) \quad \text{with } W = \mathfrak{U}_{\phi_\zeta,\arg(i\,\omega)}$$

We can then define the following geometric parametrization of all quaternionic Möbius transformations.

$$\left\{G_{\phi,\lambda,\tau} \overset{def}{=} \mathfrak{U}_{\phi,\lambda} \circ \mathfrak{D}^{\circ\tau} \circ \mathfrak{U}_{\phi,\lambda}^* \;\middle|\; \phi,\ \lambda,\ \tau \in \mathbb{R}\right\} \qquad \text{with } M = G_{\phi_\zeta,\arg(i\,\omega),\tau_\zeta}$$

$\{G_{\phi,\lambda,\tau}(\cdot)\}$ is a group isomorphic to all rotations of the Riemann sphere, that is, all rotations of $\mathbb{R}^3$, collectively denoted as $SO(3, \mathbb{R})$. Note that the subgroup $\{\mathcal{T}_\tau\}$ is not normal in $\{G_{\phi,\lambda,\tau}\}$, that is, $\{\mathcal{T}_\tau\}$ is not invariant under conjugation by an arbitrary quaternionic Möbius transformation. This can be deduced from the fact, presented in [5], that $SO(3, \mathbb{R})$ is a simple group, i.e., containing no normal subgroups other than itself and the identity transformation.

## Continuous motion

We have referred to $\{\mathcal{T}_\tau\}$ as an orbit, a term that is also applied to the subgroup $\{g^n\}_{n\in\mathbb{Z}}$ for some element $g$ of a group. The notable difference between these related conceptions of orbit is that the index $\tau$ ranges over a continuum, whereas the index $n$ is discrete.

We seek an analogous generation of $\{\mathcal{T}_\tau\}$ of the form $\exp(\tau A)$ for some transformation $A$ defined on $\hat{\mathbb{C}}$, where $\exp(\cdot)$ is an extension of the exponential function defined by its power series, in which the powers of $A$ are of the form $A^{\circ k}$, the $k$-fold composition of $A$ with itself ([1], [4], [5]). One can show that

$$\frac{d}{d\tau}\exp(\tau A)_{\tau=0} = A\exp(\tau A)_{\tau=0} = A$$

so that the conjectured generation of $\{\mathcal{T}_\tau\}$ is

$$\mathcal{T}_\tau = \exp(\tau\mathcal{T}'_0) \text{ where } \mathcal{T}'_0 = \frac{d}{d\tau}(\mathcal{T}_\tau)_{\tau=0}$$

Indeed one can show ([1], [4], [5]) that the orbit $\{\mathcal{T}_\tau\}$ is generated in this way by $\mathcal{T}'_0$, which corresponds to the initial velocity of the rotation of the Riemann sphere about the axis defined by $[u_{\zeta,\omega}]$.

Based on such derivatives of curves (not necessarily subgroups) through the identity transformation, Lie theory defines a vector space tangent to the identity transformation and equips that space with algebraic operations, notably the Lie bracket product (which in our case can be identified with the cross-product of $\mathbb{R}^3$). Quaternionic Möbius transformations are isomorphic to the group of rotations $SO(3, \mathbb{R})$, whose Lie algebra is denoted $\mathfrak{so}(3, \mathbb{R})$. This algebra can be represented by real $3\times 3$ skew-symmetric matrices ([5]). In particular, the derivative $\mathcal{T}'_0$

$$\mathcal{T}'_0 = Q\left(\cdot, i\cos(\phi_\zeta), -\sin(\phi_\zeta)\, e^{i\arg(\omega)}\right)$$

has as its counterpart in $\mathfrak{so}(3, \mathbb{R})$ the derivative of rotational matrix $[\hat{M}]$ expressed as a function of angular coordinates.

$$\frac{d}{d\tau}[\hat{M}]_{\tau=0} = \begin{pmatrix} 0 & -\cos(\phi_\zeta) & \sin(\phi_\zeta)\cos(\arg(\omega)) \\ \cos(\phi_\zeta) & 0 & \sin(\phi_\zeta)\sin(\arg(\omega)) \\ -\sin(\phi_\zeta)\cos(\arg(\omega)) & -\sin(\phi_\zeta)\sin(\arg(\omega)) & 0 \end{pmatrix}$$

Quaternionic Möbius transformations thus provide an entry to Lie theory, which is introduced more comprehensively in [1], [4], and [5].

# Glossary of Notation

$\overset{def}{=}$ — $\alpha \overset{def}{=} \beta$ means that symbol $\alpha$ is defined by expression $\beta$
(or occasionally that symbol $\beta$ is defined by expression $\alpha$)

| | |
|---|---|
| $<\cdot,\cdot>$ | real-valued inner product defined on $\mathbb{H}$ as $<h_1, h_2> \overset{def}{=} \mathcal{Re}(h_1\, h_2^*)$ |
| $[\mathcal{Vec}(h)]$ | the standard coordinates of $\mathcal{Vec}(h)$, i.e., $\{h_1, h_2, h_3\} \in \mathbb{R}^3$ |
| $[A]$ | the matrix, with respect to standard coordinates, of the linear transformation $A$ on $\mathcal{Vec}(\mathbb{H})$ or $\mathbb{R}^3$ |
| $h^*$ | the quaternion conjugate of $h \in \mathbb{H}$: $h^* \overset{def}{=} \mathcal{Re}(h) - \mathcal{Vec}(h)$ |
| $\lVert h \rVert$ | the norm of $h \in \mathbb{H}$: $\lVert h \rVert^2 \overset{def}{=} <h, h>$ |
| $\hat{\mathbb{C}}$ | the extended complex plane, $\mathbb{C} \cup \{\infty\}$ |
| $C_q(\cdot)$ | $\Gamma_q$ restricted to $\mathcal{Vec}(\mathbb{H})$, which is a rotation of $\mathcal{Vec}(\mathbb{H})$ |
| $D(\cdot)$ | $Q\left(\cdot, e^{i\,\tau_\zeta/2}, 0\right)$, a component of $M$'s polar decomposition |
| $\mathfrak{D}^{\circ\tau}(\cdot)$ | $Q\left(\cdot, e^{i\,\tau/2}, 0\right): z \mapsto e^{i\,\tau} z$, with $\tau \in \mathbb{R}$ and $D = \mathfrak{D}^{\circ\tau_\zeta}$ |
| $G_{\phi,\lambda,\tau}(\cdot)$ | $\mathfrak{U}_{\phi,\lambda} \circ \mathfrak{D}^{\circ\tau} \circ \mathfrak{U}^*_{\phi,\lambda}$, with $\phi$, $\lambda$, $\tau \in \mathbb{R}$ and $M = G_{\phi_\zeta,\arg(i\,\omega),\tau_\zeta}$ |
| $\gamma(\cdot)$ | the bijection from $\mathbb{H}$ to $\mathbb{C}^2$ given by $q_0 + q_1\, i_1 + q_2\, i_2 + q_3\, i_3 \mapsto \begin{pmatrix} q_0 + q_3\, i \\ q_2 - q_1\, i \end{pmatrix}$ |
| $\tilde{\gamma}(\cdot)$ | the isomorphism $C_q(\cdot) \mapsto Q(\cdot, \gamma(q))$ |
| $\Gamma_q(\cdot)$ | the conjugation $h \mapsto q\, h\, q^{-1}$ for quaternions $h$, $q \in \mathbb{H}$ with $q \neq 0$ |
| $\mathbb{H}$ | the non-commutative division algebra of quaternions |
| $\mathbb{H}^\times$ | the multiplicative group of non-zero quaternions, $\mathbb{H} \setminus \{0\}$ |
| $h$, $q$ | elements of the algebra $\mathbb{H}$ of quaternions, having standard coordinates denoted as $h = h_0 + h_1\, i_1 + h_2\, i_2 + h_3\, i_3$, with $h_0, h_1, h_2, h_3 \in \mathbb{R}$ |
| $\mathcal{Im}(\cdot)$ | the imaginary part of a complex value $z$, with $\mathcal{Im}(x + i\, y) = y$ |
| $M(\cdot)$ | abbreviation for Möbius transformation $Q(\cdot, \zeta, \omega)$, or more generally $\mathcal{M}(\cdot, a, b, c, d)$ |
| $\hat{M}(\cdot)$ | $\sigma^{-1} \circ M \circ \sigma$, the transformation (rotation) of $\mathbb{R}^3$ induced by $M(\cdot)$ |
| $\mathcal{M}(\cdot, a, b, c, d)$ | a Möbius transformation of the form $z \mapsto \frac{a\,z+b}{c\,z+d}$, with $z, a, b, c, d \in \mathbb{C}$ such that $a\,d - b\,c \neq 0$ |
| $Q(\cdot, \zeta, \omega)$ | a Möbius transformation of the form $z \mapsto \frac{\zeta\, z - \omega}{\overline{\omega}\, z + \overline{\zeta}}$, with $\zeta, \omega \in \mathbb{C}$ such that $\lvert\zeta\rvert^2 + \lvert\omega\rvert^2 = 1$ |
| $\mathcal{Re}(\cdot)$ | the real part of a complex value $z$, with $\mathcal{Re}(x + i\, y) = x$, or of a quaternion $h$, with $\mathcal{Re}(h_0 + h_1\, i_1 + h_2\, i_2 + h_3\, i_3) = h_0$, where $x, y, h_0, h_1, h_2, h_3 \in \mathbb{R}$ |
| $\sigma(\cdot)$ | stereographic projection from $S^2$ to $\hat{\mathbb{C}}$: $\{\eta_1, \eta_2, \eta_3\} \mapsto \frac{\eta_1 + i\,\eta_2}{1 - \eta_3}$, with $\{0, 0, 1\} \mapsto \infty$ |
| $\sigma^{-1}(\cdot)$ | the functional inverse of $\sigma(\cdot)$, i.e., $x + i\, y \mapsto \left\{\frac{2x}{x^2+y^2+1}, \frac{2y}{x^2+y^2+1}, \frac{x^2+y^2-1}{x^2+y^2+1}\right\}$ |
| $S^2$ | the Riemann sphere, $\left\{\boldsymbol{\eta} = \{\eta_1, \eta_2, \eta_3\} \in \mathbb{R}^3 \mid \eta_1^2 + \eta_2^2 + \eta_3^2 = 1\right\}$ |
| $\mathcal{T}_\tau$ | $W \circ \mathfrak{D}^{\circ\tau} \circ W^*$, with $\tau \in \mathbb{R}$, $D = \mathfrak{D}^{\circ\tau_\zeta}$, $M = \mathcal{T}_{\tau_\zeta}$ |
| $\tau_\zeta$ | $2\arctan\left(\mathcal{Re}(\zeta), \sqrt{1 - \mathcal{Re}(\zeta)^2}\right)$, with $\zeta \in \mathbb{C}$ and $\lvert\zeta\rvert \le 1$ |
| $\phi_\zeta$ | $\arctan\left(\mathcal{Im}(\zeta), \sqrt{1 - \lvert\zeta\rvert^2}\right)$, with $\zeta \in \mathbb{C}$, $\lvert\zeta\rvert \le 1$, and $\zeta \ne \pm 1$ |
| $u_{\zeta,\omega}$ | $(-\mathcal{Im}(\omega)\, i_1 + \mathcal{Re}(\omega)\, i_2 + \mathcal{Im}(\zeta)\, i_3) \Big/ \sqrt{1 - \mathcal{Re}(\zeta)^2}$, with $\zeta, \omega \in \mathbb{C}$ and $\lvert\zeta\rvert^2 + \lvert\omega\rvert^2 = 1$ |
| $v_{\zeta,\omega}$ | $\left(\mathcal{Im}(\zeta)\, (\mathcal{Im}(\omega)\, i_1 - \mathcal{Re}(\omega)\, i_2) + \lvert\omega\rvert^2\, i_3\right) \Big/ \lvert\omega\rvert \sqrt{1 - \mathcal{Re}(\zeta)^2}$, with $\zeta, \omega \in \mathbb{C}$ and $\lvert\zeta\rvert^2 + \lvert\omega\rvert^2 = 1$ |
| $w_\omega$ | $(\mathcal{Re}(\omega)\, i_1 + \mathcal{Im}(\omega)\, i_2) / \lvert\omega\rvert$, with $\omega \in \mathbb{C}$ and $0 < \lvert\omega\rvert \le 1$ |
| $\mathcal{U}(\mathbb{H})$ | $\{u \in \mathcal{Vec}(\mathbb{H}) \mid \lVert u \rVert = 1\}$, or equivalently $\left\{u \in \mathbb{H} \mid u^2 = -1\right\}$ |
| $\mathfrak{U}_{\phi,\lambda}$ | $Q\left(\cdot, \cos\left(\frac{\phi}{2}\right), -\sin\left(\frac{\phi}{2}\right) e^{i\lambda}\right)$ with $W = \mathfrak{U}_{\phi_\zeta, \arg(i\,\omega)}$ |
| $\mathcal{Vec}$ | the vector part of a quaternion $h$; $\mathcal{Vec}(h) \overset{def}{=} h - \mathcal{Re}(h) = h_1\, i_1 + h_2\, i_2 + h_3\, i_3$ |
| $\mathcal{Vec}(\mathbb{H})$ | $\{\mathcal{Vec}(h) \mid h \in \mathbb{H}\}$ |
| $W(\cdot)$ | $Q\left(\cdot, \cos(\phi_\zeta/2), -\sin(\phi_\zeta/2)\, \frac{i\,\omega}{\lvert\omega\rvert}\right)$, a component of $M$'s polar decomposition |

# Derivation of Formulas

Many of the formulas and assertions in this note are derived in the references, with a few key exceptions whose derivations we now provide.

## Matrix $[\hat{M}]$

We informally referred to $\hat{M}(\cdot)$ as the $\mathbb{R}^3$ rotation induced by $M(\cdot)$. To be explicit, we now define

$$\hat{M} \overset{def}{=} \sigma^{-1} \circ M \circ \sigma$$

We then asserted $[\hat{M}]$, the matrix of $\hat{M}$ with respect to the standard $\mathbb{R}^3$ basis, to be the following.

$$[\hat{M}] = \begin{pmatrix} \mathcal{Re}\left(\zeta^2 - \omega^2\right) & -\mathcal{Im}\left(\zeta^2 + \omega^2\right) & 2\,\mathcal{Re}\,(\zeta\,\omega) \\ \mathcal{Im}\left(\zeta^2 - \omega^2\right) & \mathcal{Re}\left(\zeta^2 + \omega^2\right) & 2\,\mathcal{Im}\,(\zeta\,\omega) \\ -2\,\mathcal{Re}\,(\zeta\,\overline{\omega}) & 2\,\mathcal{Im}\,(\zeta\,\overline{\omega}) & \lvert\zeta\rvert^2 - \lvert\omega\rvert^2 \end{pmatrix}$$

To derive this formula, we first recall the formulas for stereographic projection, $\sigma(\eta_1, \eta_2, \eta_3)$ for $\{\eta_1, \eta_2, \eta_3\} \in S^2$, and its inverse $\sigma^{-1}(z)$ for $z \in \hat{\mathbb{C}}$.

$$\sigma(\eta_1, \eta_2, \eta_3) = \begin{cases} \frac{\eta_1 + i\,\eta_2}{1-\eta_3} & \text{for } \eta_3 < 1 \\ \infty & \text{for } \eta_3 = 1 \end{cases}$$

$$\sigma^{-1}(z) = \begin{cases} \{2\,\mathcal{Re}(z), 2\,\mathcal{Im}(z), |z|^2 - 1\} / (|z|^2 + 1) & \text{for } z \in \mathbb{C} \\ \{0, 0, 1\} & \text{for } z = \infty \end{cases}$$

We now evaluate $\hat{M}(\cdot)$ on the standard basis vectors of $\mathbb{R}^3$, which yield the respective columns of $[\hat{M}]$.

$$\begin{aligned}\hat{M}(1, 0, 0) &= \sigma^{-1}(M(\sigma(1, 0, 0))) \\ &= \sigma^{-1}(M(1)) \\ &= \sigma^{-1}\left(\frac{\zeta-\omega}{\overline{\omega}+\overline{\zeta}}\right)\end{aligned}$$

Note that

$$\begin{aligned}\frac{\zeta-\omega}{\overline{\omega}+\overline{\zeta}} &= \frac{(\zeta-\omega)(\zeta+\omega)}{|\zeta+\omega|^2} \\ &= \frac{\zeta^2-\omega^2}{|\zeta+\omega|^2}\end{aligned}$$

and

$$\left|\frac{\zeta-\omega}{\overline{\omega}+\overline{\zeta}}\right| = \frac{|\zeta-\omega|}{\left|\overline{\omega}+\overline{\zeta}\right|} = \frac{|\zeta-\omega|}{|\zeta+\omega|}$$

so that

$$\begin{aligned}\left|\frac{\zeta-\omega}{\zeta+\omega}\right|^2 \pm 1 &= \frac{|\zeta-\omega|^2 \pm |\zeta+\omega|^2}{|\zeta+\omega|^2} \\ &= \frac{|\zeta|^2 + |\omega|^2 - 2\,\mathcal{Re}(\zeta\,\overline{\omega}) \pm (|\zeta|^2 + |\omega|^2 + 2\,\mathcal{Re}(\zeta\,\overline{\omega}))}{|\zeta+\omega|^2}\end{aligned}$$

Thus

$$\left|\frac{\zeta-\omega}{\zeta+\omega}\right|^2 - 1 = -4\,\frac{\mathcal{Re}(\zeta\,\overline{\omega})}{|\zeta+\omega|^2}$$

$$\left|\frac{\zeta-\omega}{\zeta+\omega}\right|^2 + 1 = 2\,\frac{|\zeta|^2 + |\omega|^2}{|\zeta+\omega|^2} = \frac{2}{|\zeta+\omega|^2} \qquad \text{since } |\zeta|^2 + |\omega|^2 = 1$$

Therefore

$$\begin{aligned}\hat{M}(1, 0, 0) &= \sigma^{-1}\left(\frac{\zeta-\omega}{\overline{\omega}+\overline{\zeta}}\right) \\ &= \begin{pmatrix} \frac{2\,\mathcal{Re}(\zeta^2-\omega^2)}{|\zeta+\omega|^2} \\ \frac{2\,\mathcal{Im}(\zeta^2-\omega^2)}{|\zeta+\omega|^2} \\ -4\,\frac{\mathcal{Re}(\zeta\,\overline{\omega})}{|\zeta+\omega|^2} \end{pmatrix} \Big/ \left(\frac{2}{|\zeta+\omega|^2}\right) \\ &= \begin{pmatrix} \mathcal{Re}(\zeta^2 - \omega^2) \\ \mathcal{Im}(\zeta^2 - \omega^2) \\ -2\,\mathcal{Re}(\zeta\,\overline{\omega}) \end{pmatrix}\end{aligned}$$

which is in accord with the purported first column of $[\hat{M}]$. Similarly we have

$$\begin{aligned}\hat{M}(0, 1, 0) &= \sigma^{-1}(M(\sigma(0, 1, 0))) \\ &= \sigma^{-1}(M(i))\end{aligned}$$

$$= \sigma^{-1}\left(\frac{\zeta\,\mathbb{i}-\omega}{\overline{\omega}\,\mathbb{i}+\overline{\zeta}}\right)$$

with

$$\frac{\zeta\,\mathbb{i}-\omega}{\overline{\omega}\,\mathbb{i}+\overline{\zeta}} = \frac{(\zeta\,\mathbb{i}-\omega)\,(\zeta-\mathbb{i}\,\omega)}{|\zeta-\mathbb{i}\,\omega|^2}$$

$$= \mathbb{i}\,\frac{(\zeta+\mathbb{i}\,\omega)\,(\zeta-\mathbb{i}\,\omega)}{|\zeta-\mathbb{i}\,\omega|^2}$$

$$= \mathbb{i}\,\frac{\zeta^2-(\mathbb{i}\,\omega)^2}{|\zeta-\mathbb{i}\,\omega|^2}$$

$$= \mathbb{i}\,\frac{\zeta^2+\omega^2}{|\zeta-\mathbb{i}\,\omega|^2}$$

$$\left|\frac{\zeta\,\mathbb{i}-\omega}{\overline{\omega}\,\mathbb{i}+\overline{\zeta}}\right| = \frac{|\zeta\,\mathbb{i}-\omega|}{\left|\overline{\omega}\,\mathbb{i}+\overline{\zeta}\right|}$$

$$= \frac{|\zeta+\mathbb{i}\,\omega|}{|\zeta-\mathbb{i}\,\omega|}$$

$$\left|\frac{\zeta+\mathbb{i}\,\omega}{\zeta-\mathbb{i}\,\omega}\right|^2 \pm 1 = \frac{|\zeta+\mathbb{i}\,\omega|^2 \pm |\zeta-\mathbb{i}\,\omega|^2}{|\zeta-\mathbb{i}\,\omega|^2}$$

$$= \frac{|\zeta|^2+|\omega|^2+2\,\mathcal{Im}(\zeta\,\overline{\omega})\pm\left(|\zeta|^2+|\omega|^2-2\,\mathcal{Im}(\zeta\,\overline{\omega})\right)}{|\zeta-\mathbb{i}\,\omega|^2}$$

$$\left|\frac{\zeta+\mathbb{i}\,\omega}{\zeta-\mathbb{i}\,\omega}\right|^2 - 1 = 4\,\frac{\mathcal{Im}(\zeta\,\overline{\omega})}{|\zeta-\mathbb{i}\,\omega|^2}$$

$$\left|\frac{\zeta+\mathbb{i}\,\omega}{\zeta-\mathbb{i}\,\omega}\right|^2 + 1 = 2\,\frac{|\zeta|^2+|\omega|^2}{|\zeta-\mathbb{i}\,\omega|^2} = \frac{2}{|\zeta-\mathbb{i}\,\omega|^2}$$

Therefore

$$\hat{M}(0,\,1,\,0) = \sigma^{-1}\left(\frac{\zeta\,\mathbb{i}-\omega}{\overline{\omega}\,\mathbb{i}+\overline{\zeta}}\right)$$

$$= \begin{pmatrix} \frac{2\,\mathcal{Re}(\mathbb{i}(\zeta^2+\omega^2))}{|\zeta-\mathbb{i}\,\omega|^2} \\ \frac{2\,\mathcal{Im}(\mathbb{i}(\zeta^2+\omega^2))}{|\zeta-\mathbb{i}\,\omega|^2} \\ 4\,\frac{\mathcal{Im}(\zeta\,\overline{\omega})}{|\zeta-\mathbb{i}\,\omega|^2} \end{pmatrix} \Big/ \left(\frac{2}{|\zeta-\mathbb{i}\,\omega|^2}\right)$$

$$= \begin{pmatrix} \mathcal{Re}(\mathbb{i}(\zeta^2+\omega^2)) \\ \mathcal{Im}(\mathbb{i}(\zeta^2+\omega^2)) \\ 2\,\mathcal{Im}(\zeta\,\overline{\omega}) \end{pmatrix}$$

$$= \begin{pmatrix} -\mathcal{Im}(\zeta^2+\omega^2) \\ \mathcal{Re}(\zeta^2+\omega^2) \\ 2\,\mathcal{Im}(\zeta\,\overline{\omega}) \end{pmatrix}$$

which is in accord with the purported second column of $[\hat{M}]$. Finally we have

$$\hat{M}(0,\,0,\,1) = \sigma^{-1}(M(\sigma(0,\,0,\,1)))$$

$$= \sigma^{-1}(M(\infty))$$

$$= \sigma^{-1}\left(\frac{\zeta}{\overline{\omega}}\right) \qquad \text{for } \omega \neq 0$$

$$= \sigma^{-1}\left(\frac{\zeta\,\omega}{|\omega|^2}\right)$$

with

$$\left|\frac{\zeta}{\omega}\right|^2 - 1 = \frac{|\zeta|^2 - |\omega|^2}{|\omega|^2}$$

$$\left|\frac{\zeta}{\omega}\right|^2 + 1 = \frac{|\zeta|^2 + |\omega|^2}{|\omega|^2} = \frac{1}{|\omega|^2}$$

so that

$$\hat{M}(0,\,0,\,1) = \sigma^{-1}\left(\frac{\zeta\,\omega}{|\omega|^2}\right)$$

$$= \begin{pmatrix} \frac{2\,\mathcal{Re}(\zeta\,\omega)}{|\omega|^2} \\ \frac{2\,\mathcal{Im}(\zeta\,\omega)}{|\omega|^2} \\ \frac{|\zeta|^2 - |\omega|^2}{|\omega|^2} \end{pmatrix} \Big/ \left(\frac{1}{|\omega|^2}\right)$$

$$= \begin{pmatrix} 2\,\mathcal{Re}(\zeta\,\omega) \\ 2\,\mathcal{Im}(\zeta\,\omega) \\ |\zeta|^2 - |\omega|^2 \end{pmatrix}$$

which coincides with the purported final column of $[\hat{M}]$.

## Matrix $[C_q]$

We alluded to $C_q(\cdot)$ as the rotation of $\mathcal{Vec}(\mathbb{H})$ induced by quaternion conjugation $\Gamma_q(\cdot)$. More precisely, $C_q(\cdot)$ is the restriction of $\Gamma_q(\cdot)$ to $\mathcal{Vec}(\mathbb{H})$. Then, as shown in [2] and [5], $C_q(\cdot)$ is a linear transformation of $\mathcal{Vec}(\mathbb{H})$. Following [1], we note that $C_q(\cdot)$ can be represented as the composition of quaternion multiplication on the left by $q$ and quaternion multiplication on the right by $q^{-1}$, both operations being linear. That is, for $h$, $p$, $q \in \mathbb{H}$ we have

$$C_q = L_q \circ R_{q^{-1}}$$

where

$$L_p(h) \overset{def}{=} p\,h$$

$$R_q(h) \overset{def}{=} h\,q$$

Note that these operations commute

$$L_p \circ R_q = R_q \circ L_p$$

These linear operations have the following matrix representations (with respect to the standard basis of $\mathbb{H}$, and with matrix multiplication as usual from the left.)

$$[L_p] = \begin{pmatrix} p_0 & -p_1 & -p_2 & -p_3 \\ p_1 & p_0 & -p_3 & p_2 \\ p_2 & p_3 & p_0 & -p_1 \\ p_3 & -p_2 & p_1 & p_0 \end{pmatrix}$$

$$[R_q] = \begin{pmatrix} q_0 & q_1 & q_2 & q_3 \\ -q_1 & q_0 & -q_3 & q_2 \\ -q_2 & q_3 & q_0 & -q_1 \\ -q_3 & -q_2 & q_1 & q_0 \end{pmatrix}$$

We now replace $L_p$ by $L_q$ and we replace $R_q$ by $R_{q^{-1}}$. In addition, we assume that $\|q\|^2 = 1$, so that $q^{-1} = q^*$. Then we have

$$[\Gamma_q] = [L_q]\,[R_{q^*}]$$

$$= \begin{pmatrix} q_0 & -q_1 & -q_2 & -q_3 \\ q_1 & q_0 & -q_3 & q_2 \\ q_2 & q_3 & q_0 & -q_1 \\ q_3 & -q_2 & q_1 & q_0 \end{pmatrix} \begin{pmatrix} q_0 & q_1 & q_2 & q_3 \\ -q_1 & q_0 & -q_3 & q_2 \\ -q_2 & q_3 & q_0 & -q_1 \\ -q_3 & -q_2 & q_1 & q_0 \end{pmatrix}$$

$$= \begin{pmatrix} \|q\|^2 & 0 & 0 & 0 \\ 0 & q_0^2 + q_1^2 - q_2^2 - q_3^2 & 2\,(q_1\,q_2 - q_0\,q_3) & 2\,(q_0\,q_2 + q_1\,q_3) \\ 0 & 2\,(q_0\,q_3 + q_1\,q_2) & q_0^2 - q_1^2 + q_2^2 - q_3^2 & 2\,(q_2\,q_3 - q_0\,q_1) \\ 0 & 2\,(q_1\,q_3 - q_0\,q_2) & 2\,(q_0\,q_1 + q_2\,q_3) & q_0^2 - q_1^2 - q_2^2 + q_3^2 \end{pmatrix}$$

$$= \begin{pmatrix} 1 & 0 & 0 & 0 \\ 0 & q_0^2 + q_1^2 - q_2^2 - q_3^2 & 2\,(q_1\,q_2 - q_0\,q_3) & 2\,(q_0\,q_2 + q_1\,q_3) \\ 0 & 2\,(q_0\,q_3 + q_1\,q_2) & q_0^2 - q_1^2 + q_2^2 - q_3^2 & 2\,(q_2\,q_3 - q_0\,q_1) \\ 0 & 2\,(q_1\,q_3 - q_0\,q_2) & 2\,(q_0\,q_1 + q_2\,q_3) & q_0^2 - q_1^2 - q_2^2 + q_3^2 \end{pmatrix}$$

Matrix $[C_q]$ is thus

$$[C_q] = \begin{pmatrix} q_0^2 + q_1^2 - q_2^2 - q_3^2 & 2\,(q_1\,q_2 - q_0\,q_3) & 2\,(q_0\,q_2 + q_1\,q_3) \\ 2\,(q_0\,q_3 + q_1\,q_2) & q_0^2 - q_1^2 + q_2^2 - q_3^2 & 2\,(q_2\,q_3 - q_0\,q_1) \\ 2\,(q_1\,q_3 - q_0\,q_2) & 2\,(q_0\,q_1 + q_2\,q_3) & q_0^2 - q_1^2 - q_2^2 + q_3^2 \end{pmatrix} \quad \text{where } \|q\| = 1$$

as claimed.

## $\tilde{\gamma} : C_q(\cdot) \mapsto Q(\cdot,\, \gamma(q))$ is an isomorphism

We defined $\gamma : \mathbb{H} \to \mathbb{C}^2$ as follows.

$$\gamma(q) \overset{def}{=} \{q_0 + i\, q_3,\, q_2 - i\, q_1\} \qquad \text{for } q = q_0 + q_1\, i_1 + q_2\, i_2 + q_3\, i_3$$

Then it is easy to check that the functional inverse of $\gamma$ is the following.

$$\gamma^{-1}\{\zeta,\, \omega\} = \mathcal{Re}(\zeta) - \mathcal{Im}(\omega)\, i_1 + \mathcal{Re}(\omega)\, i_2 + \mathcal{Im}(\zeta)\, i_3 \qquad \text{for } \zeta,\, \omega \in \mathbb{C}$$

We claimed that the following mapping of parameters

$$\{\zeta,\, \omega\} = \{\zeta_q,\, \omega_q\} = \gamma(q)$$

$$q = q_{\zeta,\omega} = \gamma^{-1}\{\zeta,\, \omega\}$$

identified matrices $[\hat{M}]$ and $[C_q]$, or more precisely, that

$$[C_q] = [Q(\cdot,\, \gamma(q))]$$

We now examine this claim.  We have

$$\zeta = \zeta_q = \gamma(q)_1 = q_0 + i\, q_3$$

$$\omega = \omega_q = \gamma(q)_2 = q_2 - i\, q_1$$

so that

$$\zeta^2 = (q_0 + i\, q_3)^2 = q_0^2 - q_3^2 + 2\, i\, q_0\, q_3$$

$$\omega^2 = (q_2 - i\, q_1)^2 = q_2^2 - q_1^2 - 2\, i\, q_1\, q_2$$

$$\zeta^2 \pm \omega^2 = q_0^2 - q_3^2 \pm \left(q_2^2 - q_1^2\right) + 2\, i (q_0\, q_3 \mp q_1\, q_2)$$

$$\zeta\,\omega = (q_0 + i\,q_3)\,(q_2 - i\,q_1)$$
$$= q_0\,q_2 + q_1\,q_3 + i(q_2\,q_3 - q_0\,q_1)$$
$$\zeta\,\overline{\omega} = (q_0 + i\,q_3)\,(q_2 + i\,q_1)$$
$$= q_0\,q_2 - q_1\,q_3 + i(q_2\,q_3 + q_0\,q_1)$$

Under this mapping $[\hat{M}]$ becomes

$$[\hat{M}] = \begin{pmatrix} \mathcal{Re}\left(\zeta^2 - \omega^2\right) & -\mathcal{Im}\left(\zeta^2 + \omega^2\right) & 2\,\mathcal{Re}\,(\zeta\,\omega) \\ \mathcal{Im}\left(\zeta^2 - \omega^2\right) & \mathcal{Re}\left(\zeta^2 + \omega^2\right) & 2\,\mathcal{Im}\,(\zeta\,\omega) \\ -2\,\mathcal{Re}\,(\zeta\,\overline{\omega}) & 2\,\mathcal{Im}\,(\zeta\,\overline{\omega}) & |\zeta|^2 - |\omega|^2 \end{pmatrix}$$
$$= \begin{pmatrix} q_0^2 + q_1^2 - q_2^2 - q_3^2 & 2\,(q_1\,q_2 - q_0\,q_3) & 2\,(q_0\,q_2 + q_1\,q_3) \\ 2\,(q_0\,q_3 + q_1\,q_2) & q_0^2 - q_1^2 + q_2^2 - q_3^2 & 2\,(q_2\,q_3 - q_0\,q_1) \\ 2\,(q_1\,q_3 - q_0\,q_2) & 2\,(q_0\,q_1 + q_2\,q_3) & q_0^2 - q_1^2 - q_2^2 + q_3^2 \end{pmatrix}$$

which is identical to the expression for $[C_q]$, as desired.

Finally, we examine the claim that $\tilde{\gamma}$ is a group isomorphism, where

$$\tilde{\gamma}(C_q(\cdot)) \overset{def}{=} Q(\cdot,\, \gamma(q)) \quad \text{for } \|q\| = 1$$

For ease of notation we continue to adopt the identification

$$\{\zeta,\, \omega\} = \{\zeta_q,\, \omega_q\} = \gamma(q) \text{ with } 1 = \|q\|^2 = |\zeta|^2 + |\omega|^2$$

so that

$$Q(\cdot,\, \zeta,\, \omega) = Q(\cdot,\, \gamma(q))$$

Now, as previously noted, $\{Q(\cdot,\, \zeta,\, \omega)\}$ is isomorphic to $SU(2, \mathbb{C})/\{\pm 1\}$, with the isomorphism sending the transformation $M(\cdot)$

$$M : z \mapsto \frac{\zeta\,z - \omega}{\overline{\omega}\,z + \overline{\zeta}}$$

to the pair of unitary transformations whose respective matrices (with respect to the standard basis for $\mathbb{C}^2$) are

$$\pm\begin{pmatrix} \zeta & -\omega \\ \overline{\omega} & \overline{\zeta} \end{pmatrix}$$

Consequently, to establish our claim it suffices to establish that $\{C_q\}$ is a group and that the following mapping is an isomorphism

$$C_{\pm q} \mapsto \left\{\pm\begin{pmatrix} \zeta & -\omega \\ \overline{\omega} & \overline{\zeta} \end{pmatrix}\right\} \quad \text{where } \{\zeta,\, \omega\} = \gamma(q)$$

We begin by noting that $\{\Gamma_q\}$ is a group, with $\Gamma_1$ equal to the identity transformation on $\mathbb{H}$.

$$\Gamma_1 : h \mapsto 1\,h\,1^{-1} = h \qquad \text{for all } h \in \mathbb{H}$$
$$(\Gamma_p \circ \Gamma_q)\,(h) = p\left(q\,h\,q^{-1}\right)p^{-1} = (p\,q)\,h\,(p\,q)^{-1} = \Gamma_{p\,q}(h) \qquad p,\, q \in \mathbb{H}^{\times}$$

Note that

$$\Gamma_{-q} = \Gamma_q$$

and as previously noted

$$[\Gamma_q] = \begin{pmatrix} 1 & 0 & 0 & 0 \\ 0 & q_0^2 + q_1^2 - q_2^2 - q_3^2 & 2(q_1 q_2 - q_0 q_3) & 2(q_0 q_2 + q_1 q_3) \\ 0 & 2(q_0 q_3 + q_1 q_2) & q_0^2 - q_1^2 + q_2^2 - q_3^2 & 2(q_2 q_3 - q_0 q_1) \\ 0 & 2(q_1 q_3 - q_0 q_2) & 2(q_0 q_1 + q_2 q_3) & q_0^2 - q_1^2 - q_2^2 + q_3^2 \end{pmatrix}$$

$$= \begin{pmatrix} 1 & 0 \\ 0 & [C_q] \end{pmatrix}$$

$$= 1 \oplus [C_q]$$

It follows that $\{C_q\}$ is a group of transformations on $\mathcal{V}ec(\mathbb{H})$, with $C_{-q} = C_q$, $C_p \circ C_q = C_{p\,q}$, and $C_1$ as the identity transformation. To establish the claim that $\tilde{\gamma}$ is an isomorphism, note first that

$$\tilde{\gamma}(C_1(\cdot)) = Q(\cdot, \zeta_1, \omega_1) = Q(\cdot, 1, 0)$$

That is, $\tilde{\gamma}$ maps the identity transformation of $\mathcal{V}ec(\mathbb{H})$ to the identity transformation of $\hat{\mathbb{C}}$. It now suffices to show that

$$\begin{pmatrix} \zeta_p & -\omega_p \\ \overline{\omega_p} & \overline{\zeta_p} \end{pmatrix} \begin{pmatrix} \zeta_q & -\omega_q \\ \overline{\omega_q} & \overline{\zeta_q} \end{pmatrix} = \begin{pmatrix} \zeta_{p\,q} & -\omega_{p\,q} \\ \overline{\omega_{p\,q}} & \overline{\zeta}_{p\,q} \end{pmatrix}$$

or equivalently, that

$$\zeta_{p\,q} = \zeta_p \zeta_q - \omega_p \overline{\omega_q}$$

$$\omega_{p\,q} = \zeta_p \omega_q + \omega_p \overline{\zeta_q}$$

We have

$$\begin{aligned} p\,q &= (p_0 + p_1 \mathbb{i}_1 + p_2 \mathbb{i}_2 + p_3 \mathbb{i}_3)(q_0 + q_1 \mathbb{i}_1 + q_2 \mathbb{i}_2 + q_3 \mathbb{i}_3) \\ &= (p_0 q_0 - p_1 q_1 - p_2 q_2 - p_3 q_3) + \\ &\quad \mathbb{i}_1 (p_0 q_1 + p_1 q_0 + p_2 q_3 - p_3 q_2) + \\ &\quad \mathbb{i}_2 (p_0 q_2 - p_1 q_3 + p_2 q_0 + p_3 q_1) + \\ &\quad \mathbb{i}_3 (p_0 q_3 + p_1 q_2 - p_2 q_1 + p_3 q_0) \end{aligned}$$

Thus

$$\zeta_{p\,q} = (p_0 q_0 - p_1 q_1 - p_2 q_2 - p_3 q_3) + \mathbb{i}(p_0 q_3 + p_1 q_2 - p_2 q_1 + p_3 q_0)$$

$$\omega_{p\,q} = (p_0 q_2 - p_1 q_3 + p_2 q_0 + p_3 q_1) - \mathbb{i}(p_0 q_1 + p_1 q_0 + p_2 q_3 - p_3 q_2)$$

and

$$\begin{aligned} &\zeta_p \zeta_q - \omega_p \overline{\omega_q} \\ &\quad = (p_0 + \mathbb{i} p_3)(q_0 + \mathbb{i} q_3) - (p_2 - \mathbb{i} p_1)\overline{(q_2 - \mathbb{i} q_1)} \\ &\quad = (p_0 + \mathbb{i} p_3)(q_0 + \mathbb{i} q_3) - (p_2 - \mathbb{i} p_1)(q_2 + \mathbb{i} q_1) \\ &\quad = (p_0 q_0 - p_3 q_3) + \mathbb{i}(p_0 q_3 + p_3 q_0) - \\ &\qquad (p_2 q_2 + p_1 q_1) - \mathbb{i}(-p_1 q_2 + p_2 q_1) \\ &\quad = (p_0 q_0 - p_1 q_1 - p_2 q_2 - p_3 q_3) + \mathbb{i}(p_0 q_3 + p_1 q_2 - p_2 q_1 + p_3 q_0) \\ &\quad = \zeta_{p\,q} \end{aligned}$$

$$\begin{aligned} &\zeta_p \omega_q + \omega_p \overline{\zeta_q} \\ &\quad = (p_0 + \mathbb{i} p_3)(q_2 - \mathbb{i} q_1) + (p_2 - \mathbb{i} p_1)\overline{(q_0 + \mathbb{i} q_3)} \end{aligned}$$

$$= (p_0 + i\, p_3)(q_2 - i\, q_1) + (p_2 - i\, p_1)(q_0 - i\, q_3)$$
$$= (p_0 q_2 + p_3 q_1) + i\,(p_3 q_2 - p_0 q_1) + (p_2 q_0 - p_1 q_3) + i\,(-p_1 q_0 - p_2 q_3)$$
$$= (p_0 q_2 - p_1 q_3 + p_2 q_0 + p_3 q_1) - i\,(p_0 q_1 + p_1 q_0 + p_2 q_3 - p_3 q_2)$$
$$= \omega_{p\,q}$$

This completes our proof of the claim that $\tilde{\gamma}$ is a group isomorphism.

## $\{\zeta, \omega\}$ as functions of $\{\tau_\zeta, \phi_\zeta, \arg(\omega)\}$

We now establish the following expressions of $\zeta$ and $\omega$ as functions of the angular arguments $\{\tau_\zeta, \phi_\zeta, \arg(\omega)\}$.

$$\zeta = \cos(\tau_\zeta/2) + i \sin(\tau_\zeta/2)\cos(\phi_\zeta)$$
$$\omega = \sin(\tau_\zeta/2)\sin(\phi_\zeta)\, e^{i \arg(\omega)}$$

We have

$$\frac{\tau_\zeta}{2} = \arctan\left(\mathcal{Re}(\zeta),\ \sqrt{1 - \mathcal{Re}(\zeta)^2}\right)$$

$$\phi_\zeta = \arctan\left(\mathcal{Im}(\zeta),\ \sqrt{1 - |\zeta|^2}\right)$$
$$= \arctan(\mathcal{Im}(\zeta),\ |\omega|)$$
$$= \arctan\left(\frac{\mathcal{Im}(\zeta)}{\sqrt{\mathcal{Im}(\zeta)^2 + 1 - |\zeta|^2}},\ \frac{|\omega|}{\sqrt{\mathcal{Im}(\zeta)^2 + 1 - |\zeta|^2}}\right)$$
$$= \arctan\left(\frac{\mathcal{Im}(\zeta)}{\sqrt{1 - \mathcal{Re}(\zeta)^2}},\ \frac{|\omega|}{\sqrt{1 - \mathcal{Re}(\zeta)^2}}\right)$$

so that

$$\begin{pmatrix} \cos(\tau_\zeta/2) \\ \sin(\tau_\zeta/2) \end{pmatrix} = \begin{pmatrix} \mathcal{Re}(\zeta) \\ \sqrt{1 - \mathcal{Re}(\zeta)^2} \end{pmatrix}$$
$$\begin{pmatrix} \cos(\phi_\zeta) \\ \sin(\phi_\zeta) \end{pmatrix} = \begin{pmatrix} \mathcal{Im}(\zeta) \\ |\omega| \end{pmatrix} \Big/ \sqrt{1 - \mathcal{Re}(\zeta)^2}$$

Therefore

$$\cos(\tau_\zeta/2) + i \sin(\tau_\zeta/2)\cos(\phi_\zeta)$$
$$= \mathcal{Re}(\zeta) + i \sqrt{1 - \mathcal{Re}(\zeta)^2}\ \mathcal{Im}(\zeta) \Big/ \sqrt{1 - \mathcal{Re}(\zeta)^2}$$
$$= \mathcal{Re}(\zeta) + i\, \mathcal{Im}(\zeta) = \zeta$$

$$\sin(\tau_\zeta/2)\sin(\phi_\zeta)\, e^{i \arg(\omega)}$$
$$= \sqrt{1 - \mathcal{Re}(\zeta)^2}\ |\omega|\, e^{i \arg(\omega)} \Big/ \sqrt{1 - \mathcal{Re}(\zeta)^2}$$
$$= |\omega|\, e^{i \arg(\omega)} = \omega$$

as claimed.

## $\frac{d}{d\tau}[\hat{M}]_{\tau=0}$

We examine the claim that

$$\frac{d}{d\tau}[\hat{M}]_{\tau=0} = \begin{pmatrix} 0 & -\cos(\phi_\zeta) & \sin(\phi_\zeta)\cos(\arg(\omega)) \\ \cos(\phi_\zeta) & 0 & \sin(\phi_\zeta)\sin(\arg(\omega)) \\ -\sin(\phi_\zeta)\cos(\arg(\omega)) & -\sin(\phi_\zeta)\sin(\arg(\omega)) & 0 \end{pmatrix}$$

We begin with the expression of $[\hat{M}]$ as a function of $\{\zeta, \omega\}$.

$$[\hat{M}] = \begin{pmatrix} \mathcal{Re}(\zeta^2 - \omega^2) & -\mathcal{Im}(\zeta^2 + \omega^2) & 2\,\mathcal{Re}(\zeta\,\omega) \\ \mathcal{Im}(\zeta^2 - \omega^2) & \mathcal{Re}(\zeta^2 + \omega^2) & 2\,\mathcal{Im}(\zeta\,\omega) \\ -2\,\mathcal{Re}(\zeta\,\overline{\omega}) & 2\,\mathcal{Im}(\zeta\,\overline{\omega}) & |\zeta|^2 - |\omega|^2 \end{pmatrix}$$

We now transform this equality by expressing $\{\zeta, \omega\}$ as functions of $\{\tau, \phi_\zeta, \arg(\omega)\}$ with $\tau$ variable and $\{\phi_\zeta, \arg(\omega)\}$ fixed. We have

$$\zeta^2 = (\cos(\tau/2) + i\sin(\tau/2)\cos(\phi_\zeta))^2$$
$$= \cos^2(\tau/2) - \sin^2(\tau/2) + 2\,i\cos(\phi_\zeta)\cos(\tau/2)\sin(\tau/2)$$
$$= \cos(\tau) + i\cos(\phi_\zeta)\sin(\tau)$$
$$|\zeta|^2 = \cos^2(\tau/2) + \sin^2(\tau/2)\cos^2(\phi_\zeta)$$
$$\omega^2 = \left(\sin(\tau/2)\sin(\phi_\zeta)\,e^{i\arg(\omega)}\right)^2$$
$$= \sin^2(\tau/2)\sin^2(\phi_\zeta)\,e^{2\,i\arg(\omega)}$$
$$= \sin^2(\tau/2)\sin^2(\phi_\zeta)\,(\cos(2\arg(\omega)) + i\sin(2\arg(\omega)))$$
$$|\omega|^2 = \sin^2(\tau/2)\sin^2(\phi_\zeta)$$
$$\zeta\,\omega = (\cos(\tau/2) + i\sin(\tau/2)\cos(\phi_\zeta))\sin(\tau/2)\sin(\phi_\zeta)\,e^{i\arg(\omega)}$$
$$\zeta\,\overline{\omega} = (\cos(\tau/2) + i\sin(\tau/2)\cos(\phi_\zeta))\sin(\tau/2)\sin(\phi_\zeta)\,e^{-i\arg(\omega)}$$

We now differentiate these expressions with respect to $\tau$, and evaluate the derivatives at $\tau = 0$.

$$\frac{d}{d\tau}\{\zeta^2\}_{\tau=0} = \frac{d}{d\tau}\{\cos(\tau) + i\cos(\phi_\zeta)\sin(\tau)\}_{\tau=0}$$
$$= \{-\sin(\tau) + i\cos(\phi_\zeta)\,\mathrm{cost}(\tau)\}_{\tau=0}$$
$$= i\cos(\phi_\zeta)$$
$$\frac{d}{d\tau}\{|\zeta|^2\}_{\tau=0} = \frac{d}{d\tau}\left\{\cos^2\left(\tfrac{\tau}{2}\right) + \sin^2\left(\tfrac{\tau}{2}\right)\cos^2(\phi_\zeta)\right\}_{\tau=0}$$
$$= \left\{-\cos\left(\tfrac{\tau}{2}\right)\sin\left(\tfrac{\tau}{2}\right) + \sin\left(\tfrac{\tau}{2}\right)\cos\left(\tfrac{\tau}{2}\right)\cos^2(\phi_\zeta)\right\}_{\tau=0}$$
$$= 0$$
$$\frac{d}{d\tau}\{\omega^2\}_{\tau=0} = \frac{d}{d\tau}\left\{\sin^2\left(\tfrac{\tau}{2}\right)\sin^2(\phi_\zeta)\,e^{2\,i\arg(\omega)}\right\}_{\tau=0}$$
$$= \left\{\sin\left(\tfrac{\tau}{2}\right)\cos\left(\tfrac{\tau}{2}\right)\sin^2(\phi_\zeta)\,e^{2\,i\arg(\omega)}\right\}_{\tau=0}$$
$$= 0$$

$$\frac{d}{d\tau}\left\{|\omega|^2\right\}_{\tau=0} = \frac{d}{d\tau}\left\{\sin^2\left(\frac{\tau}{2}\right)\sin^2(\phi_\zeta)\right\}_{\tau=0}$$
$$= \left\{\sin\left(\frac{\tau}{2}\right)\cos\left(\frac{\tau}{2}\right)\sin^2(\phi_\zeta)\right\}_{\tau=0}$$
$$= 0$$

$$\frac{d}{d\tau}\{\zeta\,\omega\}_{\tau=0} = \frac{d}{d\tau}\left\{\left(\cos\left(\frac{\tau}{2}\right) + i\sin\left(\frac{\tau}{2}\right)\cos(\phi_\zeta)\right)\sin\left(\frac{\tau}{2}\right)\sin(\phi_\zeta)\,e^{i\arg(\omega)}\right\}_{\tau=0}$$
$$= \left\{\frac{1}{2}\left(-\sin\left(\frac{\tau}{2}\right) + i\cos\left(\frac{\tau}{2}\right)\cos(\phi_\zeta)\right)\sin\left(\frac{\tau}{2}\right)\sin(\phi_\zeta)\,e^{i\arg(\omega)}\right\}_{\tau=0} + \left\{\frac{1}{2}\left(\cos\left(\frac{\tau}{2}\right) + i\sin\left(\frac{\tau}{2}\right)\cos(\phi_\zeta)\right)\cos\left(\frac{\tau}{2}\right)\sin(\phi_\zeta)\,e^{i\arg(\omega)}\right\}_{\tau=0}$$
$$= \left\{\frac{1}{2}\,(i\cos(\phi_\zeta))\times 0\right\} + \left\{\frac{1}{2}\sin(\phi_\zeta)\,e^{i\arg(\omega)}\right\}$$
$$= \frac{1}{2}\sin(\phi_\zeta)\,e^{i\arg(\omega)}$$
$$= \frac{1}{2}\sin(\phi_\zeta)\,(\cos(\arg(\omega)) + i\sin(\arg(\omega)))$$

$$\frac{d}{d\tau}\{\zeta\,\overline{\omega}\}_{\tau=0} = \frac{d}{d\tau}\left\{\left(\cos\left(\frac{\tau}{2}\right) + i\sin\left(\frac{\tau}{2}\right)\cos(\phi_\zeta)\right)\sin\left(\frac{\tau}{2}\right)\sin(\phi_\zeta)\,e^{-i\arg(\omega)}\right\}_{\tau=0}$$
$$= \left\{\frac{1}{2}\left(-\sin\left(\frac{\tau}{2}\right) + i\cos\left(\frac{\tau}{2}\right)\cos(\phi_\zeta)\right)\sin\left(\frac{\tau}{2}\right)\sin(\phi_\zeta)\,e^{-i\arg(\omega)}\right\}_{\tau=0} + \left\{\frac{1}{2}\left(\cos\left(\frac{\tau}{2}\right) + i\sin\left(\frac{\tau}{2}\right)\cos(\phi_\zeta)\right)\cos\left(\frac{\tau}{2}\right)\sin(\phi_\zeta)\,e^{-i\arg(\omega)}\right\}_{\tau=0}$$
$$= \left\{\frac{1}{2}\,(i\cos(\phi_\zeta))\times 0\right\} + \left\{\frac{1}{2}\sin(\phi_\zeta)\,e^{-i\arg(\omega)}\right\}$$
$$= \frac{1}{2}\sin(\phi_\zeta)\,e^{-i\arg(\omega)}$$
$$= \frac{1}{2}\sin(\phi_\zeta)\,(\cos(\arg(\omega)) - i\sin(\arg(\omega)))$$

Thus

$$\frac{d}{d\tau}[\hat{M}]_{\tau=0} = \frac{d}{d\tau}\begin{pmatrix} \mathcal{Re}\left(\zeta^2 - \omega^2\right) & -\mathcal{Im}\left(\zeta^2 + \omega^2\right) & 2\,\mathcal{Re}\,(\zeta\,\omega) \\ \mathcal{Im}\left(\zeta^2 - \omega^2\right) & \mathcal{Re}\left(\zeta^2 + \omega^2\right) & 2\,\mathcal{Im}\,(\zeta\,\omega) \\ -2\,\mathcal{Re}\,(\zeta\,\overline{\omega}) & 2\,\mathcal{Im}\,(\zeta\,\overline{\omega}) & |\zeta|^2 - |\omega|^2 \end{pmatrix}_{\tau=0}$$
$$= \begin{pmatrix} 0 & -\cos(\phi_\zeta) & \sin(\phi_\zeta)\cos(\arg(\omega)) \\ \cos(\phi_\zeta) & 0 & \sin(\phi_\zeta)\sin(\arg(\omega)) \\ -\sin(\phi_\zeta)\cos(\arg(\omega)) & -\sin(\phi_\zeta)\sin(\arg(\omega)) & 0 \end{pmatrix}$$

as claimed.

## $\mathcal{T}'_0 = \frac{d}{d\tau}(\mathcal{T}_\tau)_{\tau=0}$

Recall the definition of $\mathcal{T}_\tau$

$$\mathcal{T}_\tau \overset{def}{=} W \circ \mathfrak{D}^{\circ\tau} \circ W^*$$

with the following matrix representation, $[\mathcal{T}_\tau]$

$$[\mathcal{T}_\tau] = [W] \times [\mathfrak{D}^{\circ\tau}] \times [W]^*$$

$$= \begin{pmatrix} \cos\left(\frac{\phi_\zeta}{2}\right) & \sin\left(\frac{\phi_\zeta}{2}\right)\frac{i\,\omega}{|\omega|} \\ \sin\left(\frac{\phi_\zeta}{2}\right)\frac{i\,\overline{\omega}}{|\omega|} & \cos\left(\frac{\phi_\zeta}{2}\right) \end{pmatrix} \times$$

$$\begin{pmatrix} e^{i\,\tau/2} & 0 \\ 0 & e^{-i\,\tau/2} \end{pmatrix} \times \begin{pmatrix} \cos\left(\frac{\phi_\zeta}{2}\right) & -\sin\left(\frac{\phi_\zeta}{2}\right)\frac{i\,\omega}{|\omega|} \\ -\sin\left(\frac{\phi_\zeta}{2}\right)\frac{i\,\overline{\omega}}{|\omega|} & \cos\left(\frac{\phi_\zeta}{2}\right) \end{pmatrix}$$

$$= \begin{pmatrix} \cos\left(\frac{\phi_\zeta}{2}\right) & \sin\left(\frac{\phi_\zeta}{2}\right)\frac{i\,\omega}{|\omega|} \\ \sin\left(\frac{\phi_\zeta}{2}\right)\frac{i\,\overline{\omega}}{|\omega|} & \cos\left(\frac{\phi_\zeta}{2}\right) \end{pmatrix} \times$$

$$\begin{pmatrix} \cos\left(\frac{\phi_\zeta}{2}\right)e^{i\,\tau/2} & -\sin\left(\frac{\phi_\zeta}{2}\right)\frac{i\,\omega}{|\omega|}e^{i\,\tau/2} \\ -\sin\left(\frac{\phi_\zeta}{2}\right)\frac{i\,\overline{\omega}}{|\omega|}e^{-i\,\tau/2} & \cos\left(\frac{\phi_\zeta}{2}\right)e^{-i\,\tau/2} \end{pmatrix}$$

$$= \begin{pmatrix} \cos^2\left(\frac{\phi_\zeta}{2}\right)e^{i\,\tau/2} + \sin^2\left(\frac{\phi_\zeta}{2}\right)e^{-i\,\tau/2} & \cos\left(\frac{\phi_\zeta}{2}\right)\sin\left(\frac{\phi_\zeta}{2}\right)\frac{i\,\omega}{|\omega|}\left(e^{-i\,\tau/2} - e^{i\,\tau/2}\right) \\ \cos\left(\frac{\phi_\zeta}{2}\right)\sin\left(\frac{\phi_\zeta}{2}\right)\frac{i\,\overline{\omega}}{|\omega|}\left(e^{i\,\tau/2} - e^{-i\,\tau/2}\right) & \cos^2\left(\frac{\phi_\zeta}{2}\right)e^{-i\,\tau/2} + \sin^2\left(\frac{\phi_\zeta}{2}\right)e^{i\,\tau/2} \end{pmatrix}$$

$$= \begin{pmatrix} \cos\left(\frac{\tau}{2}\right) + i\sin\left(\frac{\tau}{2}\right)\cos(\phi_\zeta) & -i\sin\left(\frac{\tau}{2}\right)\sin(\phi_\zeta)\frac{i\,\omega}{|\omega|} \\ i\sin\left(\frac{\tau}{2}\right)\sin(\phi_\zeta)\frac{i\,\overline{\omega}}{|\omega|} & \cos\left(\frac{\tau}{2}\right) - i\sin\left(\frac{\tau}{2}\right)\cos(\phi_\zeta) \end{pmatrix}$$

$$= \begin{pmatrix} \cos\left(\frac{\tau}{2}\right) + i\sin\left(\frac{\tau}{2}\right)\cos(\phi_\zeta) & \sin\left(\frac{\tau}{2}\right)\sin(\phi_\zeta)\frac{\omega}{|\omega|} \\ -\sin\left(\frac{\tau}{2}\right)\sin(\phi_\zeta)\frac{\overline{\omega}}{|\omega|} & \cos\left(\frac{\tau}{2}\right) - i\sin\left(\frac{\tau}{2}\right)\cos(\phi_\zeta) \end{pmatrix}$$

Thus

$$\mathcal{T}_\tau = Q\left(\cdot,\ \cos\left(\frac{\tau}{2}\right) + i\sin\left(\frac{\tau}{2}\right)\cos(\phi_\zeta),\ -\sin\left(\frac{\tau}{2}\right)\sin(\phi_\zeta)\frac{\omega}{|\omega|}\right)$$

$$= Q\left(\cdot,\ \cos\left(\frac{\tau}{2}\right) + i\sin\left(\frac{\tau}{2}\right)\cos(\phi_\zeta),\ -\sin\left(\frac{\tau}{2}\right)\sin(\phi_\zeta)e^{i\arg(\omega)}\right)$$

Now we differentiate each matrix element with respect to $\tau$ and then set $\tau = 0$.

$$\frac{d}{d\tau}[\mathcal{T}_\tau]_{\tau=0}$$

$$= \frac{1}{2}\begin{pmatrix} -\sin\left(\frac{\tau}{2}\right) + i\cos\left(\frac{\tau}{2}\right)\cos(\phi_\zeta) & \cos\left(\frac{\tau}{2}\right)\sin(\phi_\zeta)\frac{\omega}{|\omega|} \\ -\cos\left(\frac{\tau}{2}\right)\sin(\phi_\zeta)\frac{\overline{\omega}}{|\omega|} & -\sin\left(\frac{\tau}{2}\right) - i\cos\left(\frac{\tau}{2}\right)\cos(\phi_\zeta) \end{pmatrix}_{\tau=0}$$

$$= \frac{1}{2}\begin{pmatrix} i\cos(\phi_\zeta) & \sin(\phi_\zeta)\frac{\omega}{|\omega|} \\ -\sin(\phi_\zeta)\frac{\overline{\omega}}{|\omega|} & -i\cos(\phi_\zeta) \end{pmatrix}$$

Recalling that for $k \neq 0$ the $2\times 2$ matrices $k\,A$ and $A$ refer to the same Möbius transformation, we have

$$\mathcal{T}'_0 = Q\left(\cdot,\ i\cos(\phi_\zeta),\ -\sin(\phi_\zeta)\frac{\omega}{|\omega|}\right)$$

$$= Q\left(\cdot,\ i\cos(\phi_\zeta),\ -\sin(\phi_\zeta)e^{i\arg(\omega)}\right)$$

## ■ $\sigma^{-1} \circ \mathcal{T}'_0 \circ \sigma \neq \frac{d}{d\tau}[\hat{M}]_{\tau=0}$

To simplify calculations, one might hope that the operation of differentiation would commute with stereo-graphic conjugation, which transforms functions on $\hat{\mathbb{C}}$ to functions on $S^2$. That is, we might hope that

$$\sigma^{-1} \circ \frac{d}{d\tau} (\mathcal{T}_\tau)_{\tau=0} \circ \sigma \stackrel{?}{=} \frac{d}{d\tau} \left(\sigma^{-1} \circ \mathcal{T}_\tau \circ \sigma\right)_{\tau=0}$$

To show that this is not generally the case, consider the quaternionic Möbius transformation

$$M(z) = \mathcal{T}_\tau(z) = \frac{\cos\left(\frac{\tau}{2}\right) z + i \sin\left(\frac{\tau}{2}\right)}{i \sin\left(\frac{\tau}{2}\right) z + \cos\left(\frac{\tau}{2}\right)}$$

Note that

$$M(1) = 1$$

$$\begin{aligned} M(i) &= \frac{\cos\left(\frac{\tau}{2}\right) i + i \sin\left(\frac{\tau}{2}\right)}{i \sin\left(\frac{\tau}{2}\right) i + \cos\left(\frac{\tau}{2}\right)} \\ &= i \, \frac{\cos\left(\frac{\tau}{2}\right) + \sin\left(\frac{\tau}{2}\right)}{\cos\left(\frac{\tau}{2}\right) - \sin\left(\frac{\tau}{2}\right)} \\ &= i \, \frac{\left(\cos\left(\frac{\tau}{2}\right) + \sin\left(\frac{\tau}{2}\right)\right)\left(\cos\left(\frac{\tau}{2}\right) - \sin\left(\frac{\tau}{2}\right)\right)}{\left(\cos\left(\frac{\tau}{2}\right) - \sin\left(\frac{\tau}{2}\right)\right)^2} \\ &= i \, \frac{\cos^2\left(\frac{\tau}{2}\right) - \sin^2\left(\frac{\tau}{2}\right)}{1 - 2\cos\left(\frac{\tau}{2}\right)\sin\left(\frac{\tau}{2}\right)} \\ &= i \, \frac{\cos(\tau)}{1 - \sin(\tau)} \end{aligned}$$

$$\begin{aligned} M(\infty) &= \lim_{z \to \infty} \frac{\cos\left(\frac{\tau}{2}\right) z + i \sin\left(\frac{\tau}{2}\right)}{i \sin\left(\frac{\tau}{2}\right) z + \cos\left(\frac{\tau}{2}\right)} \\ &= \frac{\cos\left(\frac{\tau}{2}\right)}{i \sin\left(\frac{\tau}{2}\right)} \\ &= \frac{-2 i \sin\left(\frac{\tau}{2}\right)\cos\left(\frac{\tau}{2}\right)}{2 \sin^2\left(\frac{\tau}{2}\right)} \\ &= \frac{-i \sin(\tau)}{1 - \cos(\tau)} \end{aligned}$$

Therefore

$$\begin{aligned} [\hat{M}] \begin{pmatrix} 1 \\ 0 \\ 0 \end{pmatrix} &= \sigma^{-1}\left(M\left(\sigma\begin{pmatrix} 1 \\ 0 \\ 0 \end{pmatrix}\right)\right) \\ &= \sigma^{-1}(M(1)) \\ &= \sigma^{-1}(1) \\ &= \begin{pmatrix} 1 \\ 0 \\ 0 \end{pmatrix} \end{aligned}$$

$$\begin{aligned} [\hat{M}] \begin{pmatrix} 0 \\ 1 \\ 0 \end{pmatrix} &= \sigma^{-1}\left(M\left(\sigma\begin{pmatrix} 0 \\ 1 \\ 0 \end{pmatrix}\right)\right) \\ &= \sigma^{-1}(M(i)) \\ &= \sigma^{-1}\left(\frac{i \cos(\tau)}{1 - \sin(\tau)}\right) \end{aligned}$$

$$= \begin{pmatrix} 0 \\ \cos(\tau) \\ \sin(\tau) \end{pmatrix}$$

$$[\hat{M}]\begin{pmatrix} 0 \\ 0 \\ 1 \end{pmatrix} = \sigma^{-1}\left(M\left(\sigma\begin{pmatrix} 0 \\ 0 \\ 1 \end{pmatrix}\right)\right)$$

$$= \sigma^{-1}(M(\infty))$$

$$= \sigma^{-1}\left(\frac{-i\sin(\tau)}{1-\cos(\tau)}\right)$$

$$= \begin{pmatrix} 0 \\ -\sin(\tau) \\ \cos(\tau) \end{pmatrix}$$

so that

$$[\hat{M}] = \begin{pmatrix} 1 & 0 & 0 \\ 0 & \cos(\tau) & -\sin(\tau) \\ 0 & \sin(\tau) & \cos(\tau) \end{pmatrix}$$

Hence

$$\frac{d}{d\tau}[\hat{M}]_{\tau=0} = \begin{pmatrix} 0 & 0 & 0 \\ 0 & -\sin(\tau) & -\cos(\tau) \\ 0 & \cos(\tau) & -\sin(\tau) \end{pmatrix}_{\tau=0}$$

$$= \begin{pmatrix} 0 & 0 & 0 \\ 0 & 0 & -1 \\ 0 & 1 & 0 \end{pmatrix}$$

On the other hand

$$\frac{d}{d\tau}(M(z))_{\tau=0} = \mathcal{T}'_0(z) = \left(\frac{-\sin\left(\frac{\tau}{2}\right)z + i\cos\left(\frac{\tau}{2}\right)}{i\cos\left(\frac{\tau}{2}\right)z - \sin\left(\frac{\tau}{2}\right)}\right)_{\tau=0}$$

$$= \frac{i}{i\,z} = \frac{1}{z}$$

with

$$\mathcal{T}'_0(1) = 1$$

$$\mathcal{T}'_0(i) = -i$$

$$\mathcal{T}'_0(\infty) = 0$$

and consequently

$$[\hat{\mathcal{T}}'_0] = \begin{pmatrix} 1 & 0 & 0 \\ 0 & -1 & 0 \\ 0 & 0 & 0 \end{pmatrix}$$

We can represent these mappings as follows

$$\begin{array}{ccc} \mathcal{T}_\tau & \xrightarrow{\sigma^{-1}\circ\mathcal{T}_\tau\circ\sigma} & \hat{\mathcal{T}}_\tau \in \mathrm{SO}(3,\mathbb{R}) \\ \frac{d}{d\tau}\Big|_{\tau=0}\downarrow & \square & \downarrow\frac{d}{d\tau}\Big|_{\tau=0} \\ \mathcal{T}'_0 & \square & \mathfrak{so}(3,\mathbb{R}) \end{array}$$

In general, then

$$\hat{\mathcal{T}}_\tau \overset{def}{=} \left(\sigma^{-1} \circ \mathcal{T}_\tau \circ \sigma\right) \in \mathrm{SO}(3, \mathbb{R})$$

so that

$$\frac{d}{d\tau}\left(\sigma^{-1} \circ \mathcal{T}_\tau \circ \sigma\right)_{\tau=0} \in \mathfrak{so}(3, \mathbb{R})$$

whereas

$$\left(\sigma^{-1} \circ \mathcal{T}'_0 \circ \sigma\right) \notin \mathfrak{so}(3, \mathbb{R})$$